\documentclass{article}

\usepackage{amsmath,amssymb,verbatim,euscript,enumerate,theorem,ifthen}

\makeatletter

\def\bddots{\mathinner{\mkern1mu\raise\p@
    \vbox{\kern7\p@\hbox{.}}\mkern2mu
        \raise4\p@\hbox{.}\mkern2mu\raise7\p@\hbox{.}\mkern1mu}}

\makeatother

\newenvironment{proof}
{
\noindent {\it Proof.}}
{\nolinebreak{\hfill$\Box\!\!$}
\bigskip}

\newenvironment{mainthm}
{\bigskip
\noindent {\bf Main Theorem}
\it 
}

\newtheorem{thm}{Theorem}[section]
\newtheorem{conjecture}[thm]{Conjecture}
\newtheorem{corollary}[thm]{Corollary}
\newtheorem{lemma}[thm]{Lemma}
\newtheorem{proposition}[thm]{Proposition}

{\theorembodyfont{\rm}
\newtheorem{definition}[thm]{Definition}

\newtheorem{remark}[thm]{Remark}
}

\numberwithin{equation}{section}

\newcommand{\A}{\mathbb{A}}

\newcommand{\J}{\mathbb J}
\newcommand{\EuA}{\EuScript A}

\newcommand{\EuC}{\EuScript C}

\newcommand{\EuH}{\EuScript H}

\newcommand{\EuZ}{\EuScript Z}
\newcommand{\fraka}{\mathfrak a}

\newcommand{\frakj}{\mathfrak j}

\newcommand{\frako}{\mathfrak o}

\newcommand{\frakt}{\mathfrak t}

\newcommand{\R}{\mathbb R}

\DeclareMathOperator{\ad}{ad}
\DeclareMathOperator{\Aff}{Aff}
\DeclareMathOperator{\aff}{\mathfrak{aff}}
\DeclareMathOperator{\Ann}{Ann}
\DeclareMathOperator{\Aut}{Aut}

\DeclareMathOperator{\Der}{Der}

\DeclareMathOperator{\gl}{\mathfrak gl}
\DeclareMathOperator{\GL}{GL}
\DeclareMathOperator{\iDer}{iDer}

\DeclareMathOperator{\Lie}{Lie}
\DeclareMathOperator{\sDer}{sDer}

\DeclareMathOperator{\spann}{span}

\DeclareMathOperator{\tr}{tr}

\begin{document}

\title{A characterization of Cayley Hypersurface and Eastwood and
  Ezhov conjecture
}

\author{}
\date{}
\maketitle

%%%%%%%%%%%%%%%%%%%%%%%%%%%%%%%%%%%%%%%%%%%%%%%%%%%%%%%%%%%%%%%%%%%%%%%%%%

\vspace{-15mm}
\begin{center}
Yuncherl Choi, \\ 
Division of General Education, Kwangwoon University \\
447-1 Wolgye-dong Nowon-Gu Seoul 139-701, Republic of Korea \\
yuncherl@kw.ac.kr, 

\vspace{3mm}
\noindent
Hyuk Kim, \\
Department of Mathematical Sciences, Seoul National University \\
San56-1 Shinrim-dong Kwanak-gu Seoul 151-747, Korea \\
hyukkim@math.snu.ac.kr
\end{center}

\vspace{3mm}
\begin{center}
\begin{minipage}[t]{105mm}\small
Eastwood and Ezhov generalized the Cayley surface to the Cayley
hypersurface in each dimension, proved some characteristic properties
of the Cayley hypersurface and conjectured that a homogeneous
hypersurface in affine space satisfying these properties must be the
Cayley hypersurface. 
We will prove this conjecture when the domain bounded by a graph of a
function defined on $\R^n$ is also homogeneous giving a
characterization of Cayley hypersurface. 
The idea of the proof is to look at the problem of affine homogeneous
hypersurfaces as that of left symmetric algebras with a Hessian type
inner product.
This method gives a new insight and powerful algebraic tools for the
study of homogeneous affine hypersurfaces.
\end{minipage}
\end{center}
%%%%%%%%%%%%%%%%%%%%%%%%%%%%%%%%%%%%%%%%%%%%%%%%%%%%%%%%%%%%%%%%%%%%%%%%%%

\vspace{2mm}
\begin{center}
\begin{minipage}[t]{105mm}

\vspace{2mm}
% 53A15 : Affine differential geometry
% 17D25 : Lie admissible algebra
% 53B05 : Linear and affine connection
\end{minipage}
\end{center}

%%%%%%%%%%%%%%%%%%%%%%%%%%%%%%%%%%%%%%%%%%%%%%%%%%%%%%%%%%%%%%%%%%%%%%%%%%

\section{Introduction}

The Cayley surface in 3-dimensional affine space is a nondegenerate
homogeneous affine surface, which  is given by
\[ z = xy - \frac13 x^3. \]
It is affinely homogeneous, that is, there is a Lie subgroup of
the group of all affine transformations of $\R^3$ that acts
transitively on the surface,  moreover its automorphism group contains
a 2-dimensional abelian Lie subgroup which acts on the surface simply
transitively\cite{NS94}.  
On a nondegenerate affine hypersurface, the cubic form $C$ is given by
$C(X,Y,Z) = (\nabla_X h)(Y,Z)$
where $\nabla$ is the induced affine connection and $h$ is the 
second fundamental form.
If the cubic form vanishes, Pick and Berwald theorem says that the
hypersurface is an open part of a nondegenerate quadric.
The Cayley surface is the first example whose cubic form $C$ is
nonzero and parallel,  that is, $C\neq 0$ and $\nabla C = 0$ \cite{NS91}.
In fact, it is the only surface satisfying this condition up to equiaffine
congruence in $\R^3$ \cite{NS94,NP89}.
For generalizations of the Cayley surface, some authors studied the
parallel cubic form of an affine hypersurface in
$\R^{n+1}$\cite{BNS90,DV90,DVY94,Tan98}.
In \cite{DV98}, Dillen and Vrancken studied the hypersurfaces with
$\nabla K = 0$, where $K$ is the difference tensor of $\nabla$ and the
Levi-Civita connection $\widehat{\nabla}$ of $h$, to obtain a very
interesting generalization. 

On the other hand,  Eastwood and Ezhov
 generalized the Cayley surface with a somewhat different view point in \cite{EE00}. 
Considering the properties of the automorphism group of the Cayley
 surface, they generalized and constructed a \emph{Cayley
hypersurface} which is homogeneous and unique up to affine congruence
 in each dimension.  
The Cayley hypersurfaces are defined as the zero
set of the following polynomial function :  
\begin{equation}\label{CayleyEqn}
\Phi(x_1, \cdots, x_{n+1}) = \sum_{d=1}^{n+1} (-1)^d \frac1d
\sum_{i+j+\cdots + m = n+1} \overbrace{x_ix_j \cdots x_m}^d .
\end{equation}
Note that, 
in the equation (\ref{CayleyEqn}), the last
coordinate occurs only once and it is the unique term of degree 1, that is,
$ \Phi(x_1, \cdots, x_n, x_{n+1}) = F(x_1, \cdots,x_n) - x_{n+1}, $
where 
\begin{equation}\label{Cayley polynomial}
F(x_1, \cdots,x_n) = \displaystyle \sum_{d=2}^{n+1} (-1)^d \frac1d
\!\!\!\! \sum_{i+j+ \cdots +m =n+1} \!\!\!\!\!\!
\overbrace{x_i x_j \cdots x_m}^d.
\end{equation}
So the Cayley hypersurface is a graph of a polynomial function $F$, which is
defined on $\R^n$,
and this polynomial will be called \emph{Cayley polynomial} in this
paper.
We will simply call a Cayley hypersurface a hypersurface which is
affinely congruent to the graph of the polynomial function
(\ref{Cayley polynomial}). 
Eastwood and Ezhov showed that a Cayley hypersurface has
the following properties  :
\begin{enumerate}[$\;\;\;\;$(E1)] 
\item The affine automorphism group contains a transitive abelian subgroup.

\item The isotropy subgroup of affine automorphism group is
 1-dimensional.

\item The affine normals are everywhere parallel. 
\end{enumerate}
Then they made a remarkable conjecture that the above properties
actually characterize the Cayley hypersurface among nondegenerate
homogeneous hypersurfaces. 
It can be easily shown that the subgroups which appeare in (E1) and (E2)
together acts transitively on the domain bounded by the Cayley
hypersurface. 
We will prove their conjecture using this additional condition. 

The following proposition is well known in affine
differential geometry and in fact follows immediately from
\cite[p. 47]{NS94}.

\begin{proposition}\label{DOH}
Let $x_{n+1} = F(x_1,\cdots, x_n)$ be a differentiable function on
$\R^n$ and let $\Sigma$ be the graph of $F$.
Then $\xi= (0, \cdots, 0,1)$ is the affine normal field of
$\Sigma$ if and only if the absolute value of the determinant of the
Hessian of $F$, $|\det DdF| = 1$. 
\end{proposition}

In the following, we will consider the nondegenerate homogeneous
hypersurface $\Sigma$ in $\R^{n+1}$, whose affine normals are given by
$\xi = (0, \cdots, 0, 1)$ and $\Sigma$ is a graph of a function which
is defined on $\R^n$. 
We will then study simply transitive subgroup of the automorphism
group of $\Sigma$, and the induced group action on $\R^n$ which may be
called a shadow action.  
It turns out that the induced connection on $\R^n$ is the standard
affine flat connection $\nabla$ and the induced metric from the
second fundamental form $h$ is invariant
by the shadow action.
Since the shadow action on $\R^n$ is simply transitive, the
standard connection gives a left symmetric algebra({\it abbreviated to} LSA)
structure on its Lie algebra, where LSA product is simply given by $x
\cdot y = \nabla_x y$ for left invariant vector fields $x$ and $y$. 
Also the invariant metric induced from $h$ gives a so called Hessian
type inner product on the Lie algebra.
In this way, a hypersurface with simply transitive action of a Lie
group corresponds to a LSA with a Hessian type inner product.  
With this correspondence between the homogeneous hypersurfaces and
Hessian LSA's, 
we will show that the abelian filiform LSA exactly corresponds to the
Cayley hypersurface. 
Then the Eastwood and Ezhov's conjecture can be transformed to the
problem of LSA.  
In this setting, the isotropy subgroup of the hypersurface appearing
in the conjecture corresponds to the Lie algebra whose elements are
derivations of the LSA as well as infinitesimal similarities with
respect to the Hessian type inner product.
Lastly, we will prove the following theorem by showing that the
dimension of the Lie algebra of the similarity derivations is affected
by the dimension of the ideal of annihilators of the LSA which in
turn is minimal only when the LSA is filiform. 

\begin{mainthm} %
Let $\Sigma$ be a nondegenerate hypersurface given by the function on
$\R^n$. Then $\Sigma$ is a Cayley hypersurface if and only if
the followings are satisfied \nolinebreak{\rm :} 
\begin{enumerate}[\rm $\;\;\;\;$(E1)]
\item $\Sigma$ admits a transitive abelian group $\A$ of affine
motions.
\item $\Aut(\Sigma)_0$  has a 1-dimensional isotropy group.
\item Affine normals to $\Sigma$ are everywhere parallel.
\item The domain above the hypersurface is also homogeneous.
\end{enumerate}
\end{mainthm}

%%%%%%%%%%%%%%%%%%%%%%%%%%%%%%%%%%%%%%%%%%%%%%%%%%%%%%%%%%%%%%%%%%%%%%%%%%%%
\section{Homogeneous hypersurface and LSA}

Let $\Sigma$  be the graph of a function $F$ defined on $\R^n$, and we
will assume that $\Sigma$ satisfies the assumptions (E1), (E2) and (E3) of Eastwood and
Ezhov's Conjecture. 
By an affine coordinate change, we may assume that $F(0) =0$, $dF_0 = 0$
and $\xi = (0,\cdots,0,1)$ is the affine normal of $\Sigma$, which says that 
$|\det DdF \: | =  1$ from Proposition \ref{DOH}. 

Let $\Aut(\Sigma)$ be the group of all affine automorphism of
$\Sigma$ and $\Aut(\Sigma)_0$ be its identity component. 
Since $F(0) =0$, i.e., $\Sigma$ contains the origin, any
element $g \in \Aut(\Sigma)$ can be written in the matrix form :
\begin{equation}\label{baseForm}
 g = \Big( \begin{pmatrix} A & b \\  c' & t \end{pmatrix},
\begin{pmatrix}x \\ F(x) \end{pmatrix} \Big) \in \Aff(n+1),
\end{equation}
where $A \in \gl(n,\R)$, $x, b, c \in \R^n $, $t \in \R$ and $c'$ is the
transpose of a column vector $c$.
Note that $g$ moves the origin to the point
$\begin{pmatrix}x \\ F(x) \end{pmatrix} \in \Sigma$.

\begin{lemma}\label{graphgroup}
With the notation in the matrix form of {\rm(\ref{baseForm})}, the affine
automorphism $g \in \Aut(\Sigma)$ is represented as the following {\rm :}
\begin{equation}\label{protoForm}
 g = \Big( \begin{pmatrix} A & 0 \\ c' & t \end{pmatrix},
\begin{pmatrix}x\\ F(x) \end{pmatrix} \Big) \in \Aff(n+1).
\end{equation}
\end{lemma}

\begin{proof}
If $g$ is equiaffine, then $g_* \xi = \xi$ since $\xi$ is equiaffine
invariant.
Let $\frakt = tI$ be a dilation, then
\[ \Sigma' = \frakt (\Sigma) = \bigg\{ \begin{pmatrix} x \\ tF(\frac{x}{t}) \end{pmatrix}
  \in \R^{n+1} \:\bigg|\: \begin{pmatrix} x \\ F(x) \end{pmatrix} \in
  \Sigma, \; t \in \R-\{ 0 \} \bigg\} \]
is a graph of $G(x) = \displaystyle t F(\frac{x}{t})$.
Since $\displaystyle |\det DdG| = \frac{1}{|t|^n} |\det DdF| =
\frac{1}{|t|^n}$ and it is constant, 
$\;|t|^{-\frac{n}{n+2}} \xi = |t|^{-\frac{2n +2}{n+2}}
\frakt_* \xi$ are affine normals of $\Sigma'$ from \cite[p.47]{NS94}.
This says that $\frakt$ preserves the direction of affine normals
  $\xi$.
Put $g = \frakt \cdot \bar{g}$ where $\frakt = tI$ for $t^n = \det g$ and $\bar{g}$
  is equiaffine. 
Then $\bar{g}(\Sigma)$ will be also a graph of a function $\bar{F}$ on
  $\R^n$. 
Since $|\det Dd\bar{F} | = t^n$, the affine normals are equal to
  $t^{\frac{n}{n+2}}\xi$ again from \cite[p.47]{NS94}. 
So $\bar{g}_*\xi = t^{\frac{n}{n+2}}\xi$, because equiaffine map
  $\bar{g}$ preserve the affine normals.
Therefore we conclude that $g$ preserves the  the direction of affine
  normals $\xi$ and hence  $b$ in (\ref{baseForm}) must vanish.  
\end{proof}

\begin{lemma} \label{affinegroup}
Let $\xi$ be the affine normal of $\Sigma$, then any element
$g$ of $\Aut(\Sigma)_0$ leaves invariant the induced connection and
acts as a conformal map with respect to the second fundamental form
$h$.
Furthermore if $g$ is equiaffine, then $g$ acts as an isometry with
respect to $h$. 
\end{lemma}

\begin{proof}
Let $D$ be the standard connection on $\R^{n+1}$ and $\nabla$ be the
induced connection on $\Sigma$. 
Since $g \in \Aut(\Sigma)_0$ is an affine map, we have 
$g_* (D_XY) = D_{g_*X}g_*Y$ for vector fields $X$ and $Y$ on $\Sigma$. 
It follows, from Gauss formula $D_XY = \nabla_XY + h(X,Y)\xi$, that
$g_*(\nabla_{X}{Y}) + h(X,Y)g_*\xi = \nabla_{g_*X}{g_*Y} + h(g_*X, g_*Y)\xi$.
Then since $g_*\xi$ is parallel to the affine normal $\xi$, we conclude that \\
\indent $\;\qquad\quad g_*(\nabla_XY) = \nabla_{g_*X} g_*Y, \quad \quad
h(X,Y)g_*\xi = h(g_*X, g_*Y)\xi. $ \\
With the matrix form (\ref{protoForm}) in Lemma \ref{graphgroup}, we
have $h(g_*X, g_*Y) = t\: h(X,Y)$, that is, $g$ is a conformal map. 
Furthermore if $g$ is equiaffine, $t =1$ since $g_*\xi = \xi$ and
hence $g$ is an isometry. 
\end{proof}

In the following, we consider  the connected component of identity
$\Aut(\Sigma)_0$, so $t$ in (\ref{protoForm}) is considered as
positive real number.  
From the action of $g$ on a point
$\begin{pmatrix} y \\ F(y) \end{pmatrix} \in \Sigma$,
we have an equation 
\begin{equation}\label{char}
F(Ay + x) = c' y + t F(y) + F(x), \qquad x,\;y \in \R^n
\end{equation}
which is the necessary condition that the matrix in (\ref{protoForm}) is an
automorphism of $\Sigma$. 
Differentiating (\ref{char}) with respect to $y$ gives us the following
equations, 
\begin{eqnarray}
\label{DiffofF}
& dF_{Ay + x} A = c'  + t dF_y, \\
\label{HessianofF}
& DdF_{Ay + x} (A\cdot, A\cdot) = t  DdF_y(\cdot, \cdot).
\end{eqnarray}
Since $|\det DdF| = 1$, (\ref{HessianofF}) says that 
\begin{equation}\label{propofAut}
\det(A)^2 = t^n, \quad t = \det(A)^{\frac{2}{n}}.
\end{equation}

Let $\A$ be an unimodular solvable subgroup of $\Aut(\Sigma)_0$ which acts on $\Sigma$
simply transitively. 
Let $q : \Sigma \rightarrow \R^n, \; \begin{pmatrix} x \\ F(x) \end{pmatrix} \mapsto
  x$ be a projection and abusing the notation, also let $q : \Aut(\Sigma)
  \rightarrow \Aff(n)$ given by  $q\Big( \begin{pmatrix} A & 0 \\ c'
    & t \end{pmatrix}, \begin{pmatrix} x \\ F(x) \end{pmatrix} \Big) =
  (A,x)$ be a projection which is a group homomorphism.
The image of the subgroup $\A$ by $q$, $\bar{\A}=q(\A)$
is an $n$ dimensional subgroup of $\Aff(n)$ and acts on $\R^n$ simply
transitively. 
For the standard flat affine connection $\bar{\nabla}$ on $\R^n$,
$\bar{\nabla} = (q^{-1})^* \nabla$  and $\bar{\nabla}$ 
induces the flat left invariant affine connection, also denoted by
$\nabla$, on $\bar{\A}$.  
Hence it defines a complete LSA structure on the Lie algebra
$\bar{\fraka} = \textnormal{Lie } \bar{\A}$ whose product is given by
$a \cdot b = \nabla_ab$ for $a, b \in \bar{\fraka}$.
In fact, since $\nabla$ is flat and torsion free, we have
\begin{eqnarray*}
& \nabla_a \nabla_b c - \nabla_b \nabla_a c - \nabla_{[a,b]} c = a(bc)
- b(ac) - [a,b]c = 0 \\
& \nabla_a b - \nabla_b a - [a,b] = ab - bc -[a,b] = 0
\end{eqnarray*}
and we obtain the left symmetry of the product, i.e., $(a,b,c) = (b,a,c)$ where the
associator $(a,b,c) = (ab)c- a(bc)$ and the compatibility with respect
to the Lie structure, $ab - ba = [a,b]$. 
We will denote the LSA as $\EuA = (\bar{\fraka}, \cdot)$ and the left
multiplication by $a\in \EuA$ as $\lambda_a$ so that $\lambda_a (b) =
ab$. 
By identifying $\EuA$ with $\R^n$, $\lambda$ induces a
representation of $\bar{\fraka}$ into $\aff(n) = \gl(n) + \R^n$ which
maps $a$ to $(\lambda_a, a)$. 
Exponentiating this representation gives us the representation of
$\bar{\A}$ into $\Aff(n) = \GL(n) \ltimes \R^n$ so that the group 
\[ \bar{\A} = \{ (\exp \lambda_a, \; a + \displaystyle
\frac{1}{2!} a^2 + \frac{1}{3!} a^3 + \cdots) \: | \: a \in
\bar{\fraka} \} \subset \Aff(n), \]
where $a^k = (\lambda_a)^{k-1} a, \; (k = 1, 2, \cdots)$, acts on
$\R^n$ simply transitively({\it see} \cite{Kim97} for details). 
If we consider the extension $\tilde{\EuA}= \EuA \oplus \R\cdot 1$ of
$\EuA$ with the identity 1, the expression $a + \displaystyle
\frac{1}{2!} a^2 + \frac{1}{3!} a^3 + \cdots$ can be denoted simply as
$e^a - 1$ in the representation of the associated Lie group.
In the following, we will use the notation $e^a - 1 = a + \displaystyle \frac{1}{2!} a^2 +
\frac{1}{3!} a^3 + \cdots$ without further comment about the extension
of a given LSA.
Since $\bar{\A}$ acts on the whole space $\R^n$, the induced LSA
$\EuA$ must be complete. 
The following proposition about the equivalent conditions of the
completeness is well known\cite{GH86, Sag92} \nolinebreak{:}

\begin{proposition} 
Let $\EuA$ be an LSA. Then the following statements are equivalent
{\rm :}
\begin{enumerate}[\hspace{4mm}\rm (a)]
\item $\EuA$ is complete.
\item $\tr \rho_a = 0$ for all $a \in \EuA$ where $\rho_a$ is the
  right multiplication of $a \in \EuA$.
\item $\det(I + \rho_a) = 1$ for all $a \in \EuA$.
\end{enumerate}
\end{proposition}

Note that $q|_{\A}$ is an isomorphism and hence $\bar{\A}$ is
unimodular. 
Then we have that $\tr \lambda_{a} = 0$ for all $a \in
\EuA$ because the Lie algebra $\bar{\fraka}$ is unimodular, that is,
$0 = \tr \ad_a = \tr \lambda_a - \tr \rho_a$.
Therefore $\det \exp \lambda_a = 1$.
Furthermore we have the following :

\begin{proposition} \label{representation}
Let $\A$ be an unimodular solvable subgroup of $\Aut(\Sigma)_0$ which acts on $\Sigma$
simply transitively. 
For any element $g \in \A$, there exists a unique $x \in \R^n$ such that
$g = g_x$ where $g_x$ is represented by
\[ g_x = \big( \begin{pmatrix} M_x & 0 \\ dF_xM_x & 1 \end{pmatrix},
\begin{pmatrix} x \\ F(x) \end{pmatrix} \big) \]
where $M_x \in \GL(n)$ is the $n \times n$ matrix given
by the equations {\rm :} for $a \in \EuA$
\[ M_x = \exp \lambda_{a}  \;\mbox{ and } \;  x= e^{a} -1
\in \R^n. \]
Moreover $(M_x,x)$ acts on $\R^n$ as an isometry with respect to the
Hessian metric $DdF$. 
\end{proposition}

\begin{proof}
Since the projection $q$ is 1-1, $\A$ is identified with $\bar{\A}=q(\A)$
which acts on $\R^n$ simply transitively.
Then by the evaluation map at $0$, $\bar{\A}$
hence $\A$ can be identified with $\R^n$.
Therefore, for any $g \in \A$, there exists a unique element
$x \in \R^n$ such that $g = g_x$ and from (\ref{protoForm}),
$g_x$ can be written as
\[ g_x = \big( \begin{pmatrix} M_x & 0 \\ c'_x & t_x \end{pmatrix},
\begin{pmatrix} x \\ F(x) \end{pmatrix} \big), \]
where $M_x \in \GL(n), \; c_x \in \R^n$ and $t_x \in \R_+$.
Note that $M_x, c_x$ and $t_x$ depend on $x \in \R^n$ smoothly.
From the canonical representation of $\EuA$ into $\gl(\R^n) + \R^n$,
there exists $a \in \EuA$ such that
$ x = e^{a} - 1 \mbox{ and } M_x = \exp \lambda_{a}$,
where $\lambda_{a}$ is a left multiplication by $a$ in the LSA
$\EuA$.
From the equation (\ref{DiffofF}) and using also that $dF_0 =
0$, we find that $c'_x = dF_x M_x$.  
Since $\det M_x = \det \exp \lambda_a = 1$, we have $t_x = 1$ by
(\ref{propofAut}). 
Therefore, we see from (\ref{HessianofF}) that 
$(M_x, x)$ acts on $\R^n$ as an isometry with respect
to the Hessian metric $DdF$.
\end{proof}

The induced affine metric $(q^{-1})^* h$ is equal to the Hessian of the
function $F$({\it see} \cite{NS94}). 
So $DdF = (q^{-1})^* h$ defines a left invariant metric on $\bar{\A}$
and induces an inner product $H = DdF_0$ on the Lie algebra
$\bar{\fraka}$.
We will frequently identify $\bar{\A}$($\EuA$, resp.) and
$\R^n$($T_0\R^n = \R^n$, resp.) via the evaluation map at $0$(its
differential, resp.) in the following so that the left invariant
vector fields on $\bar{\A}$ becomes a vector field on $\R^n$.

\begin{proposition}
The induced inner product $H$ on $\EuA$ is of Hessian type, that is,
$H$ satisfies   
\begin{equation}\label{Hessiantype}
 H(a, bc) - H(ab, c) =  H(b, ac) - H(ba, c),
\end{equation}
for all $a,b,c \in \EuA$.
\end{proposition}

\begin{proof}
From the Codazzi equation, we have
\begin{equation}\label{codazzi}
(\bar{\nabla}_a DdF_x)(b, c) = (\bar{\nabla}_b DdF_x)(a, c)
\end{equation}
for all left invariant vector fields $a,b,c \in \EuA$ and $x \in \R^n$
where $\bar{\nabla}$ is the standard flat affine connection.
The calculation at $x= 0$ gives us
\begin{eqnarray*}
\lefteqn{(\bar{\nabla}_a DdF)\big|_{x=0}(b, c)} \\ 
&=& \bar{\nabla}_a (DdF_x(b,c))\big|_{x=0} -
DdF_0(\bar{\nabla}_a b\big|_{x=0}, c) - DdF_0(b, \bar{\nabla}_a c\big|_{x=0}) \\
&=& - DdF_0(\bar{\nabla}_a b, c) - DdF_0(b, \bar{\nabla}_a c)
\end{eqnarray*}
since $b,c$ are left invariant and $DdF_x(b_x,c_x) = DdF_x(M_xb_0, M_x c_0) =
DdF_0(b_0,c_0)$ holds from (\ref{HessianofF}) and hence constant.     
Now from (\ref{codazzi}), we obtain $\displaystyle  H(ab, c) + H(b,
ac) =  H(ba, c) + H(a,bc)$.  
\end{proof}

\begin{remark} \label{remark}
The Hessian type inner product has been used by Vinberg in \cite{Vin63}.
He introduced the inner product on an LSA from the Lie algebra homomorphism
given by the trace form of left multiplication. 
Then he showed that the set of clans, which are  LSA's with the
positive definite Hessian type inner product such that the left
multiplication of any element has only real eigenvalues, is in
one-to-one correspondence with the set of the homogeneous convex
domains. 
The terminology of Hessian type was introduced by Hirohiko Shima. 
He introduced the Hessian type inner product for the study of
homogeneous Hessian manifolds\cite{Shi80}.
Using this inner product, he developed the Hessian algebra
which is the LSA with the Hessian type inner product.
\end{remark}

We have just shown that the homogeneous affine hypersurface
which is a graph of a function $F : \R^n \rightarrow \R$ with $|\det
DdF| = 1$ and the automorphism group contains a solvable unimodular
simply transitive subgroup gives a complete LSA with the Hessian type
inner product. 
We will consider the converse of this in the following.
We start from a complete LSA $\EuA$ with a Hessian
type inner product $H$ with $|\det H | = 1$.
(Abusing the notation, we
denote by $H$ in this paper both the inner product and its associated
symmetric matrix with respect to the standard basis on $\R^n$.) 
Let $\bar{\fraka}$ be the associated Lie
algebra of $\EuA$.  
Define a map $\phi : \bar{\fraka} 
\rightarrow \aff (n+1)$ by $\phi(a) = \big( \begin{pmatrix} \lambda_a
  & 0 \\ a'H & 0 \end{pmatrix} \begin{pmatrix} a \\ 0 \end{pmatrix}
\big)$,
where $a'$ is the transpose of $a \in \bar{\fraka}$.
Then $\phi$ is the Lie algebra homomorphism since, for any $a,b \in \bar{\fraka}$, 
\begin{eqnarray*}
[\phi(a), \phi(b)] &=& \phi(a)\phi(b) - \phi(b)\phi(a) \\
&=& \big( \begin{pmatrix} \lambda_a\lambda_b - \lambda_b\lambda_a
  & 0 \\ a'H\lambda_b - b'H\lambda_a  & 0 \end{pmatrix} \begin{pmatrix} 
  \lambda_a b - \lambda_b a \\ a'Hb - b'Ha \end{pmatrix} \big) \\
&=& \big( \begin{pmatrix} [\lambda_a, \lambda_b]
  & 0 \\ (ab)'H - (ba)'H & 0 \end{pmatrix} \begin{pmatrix} 
  ab - ba \\ 0 \end{pmatrix} \big) \\
&=& \big( \begin{pmatrix} \lambda_{[a,b]}
  & 0 \\ [a,b]'H & 0 \end{pmatrix} \begin{pmatrix}
  [a, b] \\ 0 \end{pmatrix} \big) \\
&=& \phi([a,b]),
\end{eqnarray*}
where the third equality follows from (\ref{Hessiantype}) and
symmetry of the Hessian inner product $H$, the forth equality
is the LSA condition and the compatibility of LSA product and Lie
product.   
Since $\fraka = \phi(\bar{\fraka}) $ is a Lie subalgebra of
$\aff(n+1)$, we have a corresponding Lie subgroup $\A = \exp \fraka 
\subset \Aff(n+1)$ whose elements are given by
\begin{equation}\label{gpelement}
 \big( \begin{pmatrix} e^{\lambda_a}
  & 0 \\ a'H(I + \frac{\lambda_a}{2!} + \frac{(\lambda_a)^2}{3!} +
  \cdots) & 1 \end{pmatrix} \begin{pmatrix} e^a -1\\
  a'H(\frac{a}{2!} + \frac{a^2}{3!} +
  \cdots) \end{pmatrix} \big), 
\end{equation}
for $a \in \bar{\fraka}$. 
Note that the orbit space of the Lie
group $\bar{\A} = \{  (e^{\lambda_a}, e^a -1) \in \aff(n)  \: | \: a
\in \EuA \}$ at the origin is the whole space $\R^n$ because $\EuA$ is
complete.
Put $x = e^a -1 \in \R^n$ and let $\displaystyle F(x) = F(e^a - 1) =
a'H(\frac{a}{2!} + \frac{a^2}{3!} + \cdots)$ be the function from
$\R^n$ to $\R$.  
Then the function $F$ defines a homogeneous
hypersurface $\Sigma = \{ (x,F(x))\in \R^{n+1} \: | \: x \in \R^n \}$
containing the origin of $\R^{n+1}$ since $F(0) =
F(e^0-1) = 0$ and  
$\A$ acts simply transitively on $\Sigma$ since $\Sigma$ is the
$\A$-orbit of the origin.  
For $a, b \in \EuA$ and $u \in \R$, the group multiplication in $\A$
shows that
\begin{eqnarray}
\lefteqn{F((e^a-1) + e^{\lambda_a}(e^{ub} -1)) = a'H(\frac{a}{2!} + \frac{a^2}{3!} +
  \cdots) }  \nonumber \\
& \qquad\qquad\qquad\qquad &   + a'H(I + \frac{\lambda_a}{2!} + \frac{(\lambda_a)^2}{3!} + 
  \cdots)(e^{ub} -1) \label{fordF}\\ 
&\qquad\qquad\qquad\qquad & + (ub)'H(\frac{ub}{2!} + \frac{(ub)^2 }{3!} +
  \cdots) .\nonumber
\end{eqnarray}
By differentiating (\ref{fordF}) with respect to $u$ at $u=0$, we obtain
\begin{equation}\label{dF}
dF(e^a -1) e^{\lambda_a} = a'H(I + \frac{\lambda_a}{2!} +
\frac{(\lambda_a)^2}{3!} + \cdots). 
\end{equation}
Therefore, if we put $M_x = e^{\lambda_a}$ where
$x = e^a -1$, then the elements of $\A$ in 
(\ref{gpelement}) is equal to the matrix in Proposition
\ref{representation}.
Differentiating (\ref{dF}) at the origin, we have
\[ DdF_0 = H. \]
Now  the group multiplication in $\A$ gives us : for $a, b \in \EuA$
and $u \in \R$,  
\begin{eqnarray}%\label{forDdF}
\lefteqn{dF((e^a-1) + e^{\lambda_a}(e^{ub} -1)) e^{\lambda_a}
  e^{\lambda_{ub}}} \nonumber \\
&=& dF(e^a-1) e^{\lambda_a} e^{\lambda_{ub}} + dF(e^{ub} -1)
e^{\lambda_{ub}}. \label{forDdF} 
\end{eqnarray}  
Once more, differentiating (\ref{forDdF}) with respect to $u$ at
$u=0$ gives 
\[ DdF_x(M_x, M_x) = DdF(e^a -1) (e^{\lambda_a}, e^{\lambda_a})  =
 DdF(0)( \cdot, \cdot) = DdF_0 (\cdot, \cdot). \]
Since $\det M_x =1$,
we have $\det DdF_x = \det DdF_0$  and hence $|\det DdF_x | = | \det H
| = 1$ for
all $x \in \R^n$. 
Therefore the affine normals of the hypersurface $\Sigma$ are equal to
$\xi = (0, \cdots, 0,1)$, so they are parallel.
Now we can summarize as follows. 
   
\begin{thm} \label{one-to-one}
There is a one-to-one correspondence between the set of the graph of
$F :\R^n \rightarrow \R$ with $|\det DdF| = 1$ on which  a unimodular Lie
subgroup $\A$ of affine transformations acts simply transitively and
the set of the complete LSA $\EuA$ with a nondegenerate Hessian
type inner product $H$ with $|\det H| = 1$.
\end{thm}

%%%%%%%%%%%%%%%%%%%%%%%%%%%%%%%%%%%%%%%%%%%%%%%%%%%%%%%%%%%%%%%%%%%%%%%%%%
\section {The transitive action of abelian subgroup}

With the first assumption (E1) of Eastwood and Ezhov conjecture, we will
assume that the Lie group $\A$ is a transitive abelian subgroup of
$\Aut(\Sigma)_0$ in this section.   
Let's denote by $\J$ the isotropy subgroup of $\Aut(\Sigma)_0$ at the
origin $0 \in \Sigma$.

\begin{proposition}\label{simply}
Let $\Sigma$ be the graph of a function $F: \R^{n}
\rightarrow \R$ with $|\det DdF| = 1$ and let $\A$ be an abelian
subgroup of $\Aut(\Sigma)_0$ which acts on $\Sigma$ transitively. Then
the action of $\A$ on $\Sigma$ is simply transitive.
\end{proposition}

\begin{proof}
Suppose that the action of $\A$ is not simply transitive.
Then $\A$ contains a nonzero element $j$ of $\J$ and by the commuting
property, we have that
$ jg\cdot 0 = gj\cdot 0 = g\cdot 0$ for all $g \in \A$. 
It follows that $g \cdot 0$ is fixed under $j$ for any $g \in \A$.
Since $\A$ acts transitively on $\Sigma$,
each $j \in \A \cap \J$ fixes every point of $\Sigma$.
Hence we see that $\Sigma$ is a hyperplane through $0$ in $\R^{n+1}$
because the fixed point set of a linear map $j$ is a subspace.
But this contradicts the nondegeneracy of the hypersurface $\Sigma$.
\end{proof}

For the abelian LSA $\EuA$, the condition (\ref{Hessiantype})
of Hessian type inner product is reduced to the following
\begin{equation}\label{HessiantypeofAbel}
H(a, bc) = H(b,ac) , \quad a,b,c \in \EuA.
\end{equation}

%By abusing the notation, we write the inner product $H(a,b)= \langle
%a, b \rangle$ as $a'Hb$ using the corresponding matrixes with respect
%to the standard basis of $\R^n$.

\begin{lemma}\label{gpelementofAbel}
Let $\Sigma$ be the graph of a function $F: \R^n \rightarrow \R$ with
$|\det DdF| = 1$ and let $\A$ be an abelian 
subgroup of $\Aut(\Sigma)_0$ which acts on $\Sigma$ simply transitively.
Then any element of $\A$ can be written as 
\[  g_x = \big( \begin{pmatrix} I + \lambda_x
  & 0 \\ x'H & 1 \end{pmatrix} \begin{pmatrix} x \\ F(x) \end{pmatrix}
  \big),\]
for some $x \in \R^n$ and $\lambda : \R^n \rightarrow \gl(\R^n)$ is a
 left multiplication map coming from the associated abelian LSA $\EuA$.
\end{lemma}

\begin{proof}
For the abelian LSA $\EuA$, we have $(\lambda_a)^n =
\lambda_{a^n}$ for all $a \in \EuA$ and $n = 2, 3, \cdots$.
Therefore,  
\begin{eqnarray*}
e^{\lambda_a} &=& I + \lambda_a + \frac{1}{2!} (\lambda_a)^2 +
\frac{1}{3!} (\lambda_a)^3 + \cdots\\ 
&=& I + \lambda_a + \frac{1}{2!} \lambda_{a^2} +
\frac{1}{3!} \lambda_{a^3} + \cdots\\ 
&=& I + \lambda_{a + \frac{1}{2!} a^2 + \frac{1}{3!} a^3 + \cdots}\\
&=& I + \lambda_x.
\end{eqnarray*}
Then from the equation (\ref{dF}) and (\ref{HessiantypeofAbel}), we have
\begin{eqnarray*}
dF_x(I + \lambda_x) &=& a'H(I + \frac{\lambda_a}{2!} +
\frac{(\lambda_a)^2}{3!} + \cdots) \\
&=& (a + \frac{1}{2!} a^2 + \frac{1}{3!} a^3 + \cdots)'H \\ 
&=& (e^a - 1)'H = x'H.
\end{eqnarray*}
Then from Proposition \ref{representation}, the elements of $\A$ can be
written in the form given above.
\end{proof} 
 
Lemma \ref{gpelementofAbel} shows that from the simply transitive abelian Lie
subgroup we can obtain the structure of the associated abelian LSA
and the Hessian type inner product
directly, and vice versa.

\begin{thm}\label{polynom}
Let $\Sigma$ be the graph of a function $F: \R^{n}
\rightarrow \R$ with $|\det DdF| = 1$. 
Suppose that
$\Aut(\Sigma)_0$ contains an abelian
subgroup which acts on $\Sigma$ simply transitively.
Then we have the followings {\rm :}
\begin{enumerate}[\hspace{3mm}\rm (a)]
\item The function $F$ is a polynomial which is given by 
\begin{equation}\label{poly}
F(x) = x'H(\frac{1}{2}I - \frac{1}{3}\lambda_x + \frac{1}{4}(\lambda_x)^2
- \cdots)x, 
\end{equation}
where $H = DdF_0$ and $\lambda$ is the left multiplication map obtained from the
associated complete abelian LSA.
\item An abelian subgroup of $\Aut(\Sigma)_0$ which acts simply transitively
on $\Sigma$ is unique. Therefore $\A$ is a normal subgroup of $\Aut(\Sigma)_0$.
\end{enumerate}
\end{thm}

\begin{proof}
(a) Let $\A$ be the simply transitive abelian subgroup of $\Aut(\Sigma)_0$ and
let $\EuA$ be the associated complete abelian LSA.  
Let $\lambda$ be the left multiplication map obtained from the LSA $\EuA$ and let
$H = DdF_0$ be the Hessian type inner product on $\EuA$. Then $H
\lambda_x = (\lambda_x)'H$ for all $x \in \R^n$ from
(\ref{HessiantypeofAbel}).  
Define the polynomial $G : \R^n \rightarrow \R$ as the right hand side of (\ref{poly}). 
If we differentiate $G(x)$ by using $d((\lambda_x)^n)(v) = n
(\lambda_x)^{n-1}\lambda_v$, then  
\begin{eqnarray*} 
dG_x &=& H(\frac{1}{2}I - \frac{1}{3}\lambda_x +
\frac{1}{4}(\lambda_x)^2 - \cdots)x \\ 
&& \qquad + x'H(\frac{1}{2}I - \frac{2}{3}\lambda_x + \frac{3}{4}(\lambda_x)^2
- \cdots)\\ 
&=&x'H(I - \lambda_x + (\lambda_x)^2 - (\lambda_x)^3 + \cdots) \\
&=& x'H(I+\lambda_x)^{-1}.
\end{eqnarray*} 
where the second equality follows from the associativity of $H$.
On the other hand, from the proof of Lemma \ref{gpelementofAbel}, 
$dF_x = x'H(I+\lambda_x)^{-1}$.
Hence we have $F(x) = G(x) + c$ for some constant $c$. 
But $G(0) = 0 = F(0)$,  so we conclude that they are equal.

(b) Using the associativity of $H$, $dF_x = x'((I+\lambda_x)^{-1})'
H$. Then, we have
\begin{equation}\label{final}
((I+\lambda_x)^{-1} x)' = dF_x H^{-1}.
\end{equation}
Because the right hand side of (\ref{final}) depends only on the given
function, we can show that the LSA structure is uniquely determined as
follows : 
If we assume that there exists a Lie algebra homomorphism
$\tilde{\lambda} : \R^n \rightarrow \gl(n)$ which gives another
associated abelian LSA structure, so that
$(I + \tilde{\lambda}_x)^{-1} x = (I + \lambda_x)^{-1} x$ for all $x
\in \R^n$ by (\ref{final}), then we have
\begin{equation}\label{differ}
[(\lambda_x - \tilde{\lambda}_x) - ((\lambda_x)^2 -
(\tilde{\lambda}_x)^2) + ((\lambda_x)^3 - (\tilde{\lambda}_x)^3) - 
\cdots ] x = 0, 
\end{equation}
for all  $x \in \R^n$. 
If we differentiate (\ref{differ}), then we can easily derive that  
\begin{equation}\label{differ2}
2(\lambda_x - \tilde{\lambda}_x) - 3 ((\lambda_x)^2 -
(\tilde{\lambda}_x)^2) + 4((\lambda_x)^3 - (\tilde{\lambda}_x)^3) - 
\cdots  = 0,
\end{equation}
for all $x \in \R^n$.
By differentiating (\ref{differ2}) at $x= 0$, we obtain $2(\lambda -
\tilde{\lambda}) = 0$, that is, $\lambda =
\tilde{\lambda}$.

Now (a) and Lemma \ref{gpelementofAbel} shows that the abelian subgroup $\A$
is also uniquely determined from the LSA structure and the Hessian
type inner product. 
\end{proof}

Any element $j \in \tilde{\J}$, the isotropy subgroup at the origin $0 \in
 \Sigma$, must be linear. So it is given by
 $ j = \begin{pmatrix} A & 0 \\ c'   & t \end{pmatrix} \in
 \GL(n+1). $ 
Since $dF_0 = 0$, the tangent space of $\Sigma$ at $0$ is $\R^n$ and
 it is invariant by the action of $j_* = j$.
Hence we have the following :

\begin{lemma}\label{elem-isot}
Any element $j \in \tilde{\J}$ is represented as the following form
\[ j = \begin{pmatrix} A & 0 \\ 0 & t \end{pmatrix}, \]
where $A \in \GL(n)$ and $t = (\det A)^{\frac{2}{n}}$.
\end{lemma}

\begin{proposition} \label{TFAE}
Let $\tilde{\J}$ be the isotropy subgroup of $\Aut(\Sigma)$ at $0$. Then
the followings are equivalent$\;${\rm :}
\begin{enumerate}[\hspace{5mm}\rm (a)]
\item $g \in \tilde{\J}$.
\item $g \in \GL(n+1)$ normalizes $\A$.
\item $g = \begin{pmatrix} A & 0 \\ 0 & t \end{pmatrix}\in \GL(n+1)$ and 
$F(Ax) = tF(x)$.
\item $g = \begin{pmatrix} A & 0 \\ 0 & t \end{pmatrix}\in
  \GL(n+1)$ and $A$ is an automorphism of the associated complete
  abelian LSA $\EuA$ and $A'HA = tH$.
\end{enumerate}
\end{proposition}

\begin{proof}
\underline{(a) $\Rightarrow$ (b) } : Since $\tilde{\J}$ is the isotropy
subgroup at 0, $g$ is contained in $\GL(n+1)$. From Theorem 
\ref{polynom}, $\A$ is the unique abelian normal subgroup of
$\Aut(\Sigma)_0$. Then $g\A g^{-1} \subset \Aut(\Sigma)_0$ must be
equal to $\A$. Hence $g$ normalizes $\A$.

\noindent
\underline{(b) $\Rightarrow$ (a),(c) and (d)} :
Let $g \in \GL(n+1)$ be a normalizer of the abelian subgroup $\A$, then for
any $h \in \A$ there exists $\tilde{h} \in \A$ such that $g\cdot (h\cdot 0) =
\tilde{h}\cdot (g \cdot 0) = \tilde{h} \cdot 0$. $\;$ 
Since $\A$ acts on $\Sigma$ transitively,  this shows that $g$ belongs to
$\Aut(\Sigma)$, that is, $g \in \tilde{\J}$. 
Therefore from Lemma 
\ref{elem-isot}, $g = \begin{pmatrix} A & 0 \\ 0 & t \end{pmatrix}$
where $A \in \GL(n)$ and $t = (\det A)^{\frac{2}{n}}$. 
Then for any element $g_x = \big( \begin{pmatrix} M_x
  & 0 \\ x' H & 1 \end{pmatrix} \begin{pmatrix} x \\ F(x) \end{pmatrix}
  \big) \in \A$, we have  
\begin{equation}\label{conju}
gg_xg^{-1} = \big( \begin{pmatrix} AM_xA^{-1}
  & 0 \\ tx' H A^{-1} & 1 \end{pmatrix} \begin{pmatrix} Ax \\ t F(x)
  \end{pmatrix} \big) \in \A.
\end{equation} 
Hence we have $A\lambda_x A^{-1} = \lambda_{Ax}$, that is, $A \in
  \Aut(\EuA)$ and $F(Ax) = tF(x)$. 
Moreover we have $tx'H A^{-1} = (Ax)'H = x'A'H$, and hence $tH =
  A'HA$. 

\noindent
\underline{(c) $\Rightarrow$ (a)} :
Since $g \cdot \begin{pmatrix} x \\  F(x) \end{pmatrix} =
\begin{pmatrix} Ax \\ t F(x) \end{pmatrix} \in \Sigma$, $g \in \tilde{\J}
\subset \Aut(\Sigma)$.

\noindent
\underline{(d) $\Rightarrow$ (b) and (c)} :
From the conditions of (d), we have
\begin{eqnarray*}
& AM_x A^{-1} = I + A\lambda_xA^{-1} = I + \lambda_{Ax} = M_{Ax} \\
& tx'HA^{-1} = x'A'H = (Ax)' H. 
\end{eqnarray*}
Since $A$ is an automorphism of
the LSA $\EuA$,
\begin{eqnarray*}
F(Ax) &=& x'A'H(\frac{1}{2}I - \frac{1}{3}A \lambda_x A^{-1} +
\frac{1}{4} A (\lambda_x)^2 A^{-1} - \cdots) Ax \\
&=& x'A'HA(\frac{1}{2}I - \frac{1}{3} \lambda_x  +
\frac{1}{4} (\lambda_x)^2  - \cdots) x \\
&=& tF(x).
\end{eqnarray*}
Using the above, (\ref{conju}) shows that $g$ normalizes the abelian
subgroup $\A$. 
\end{proof}

From Proposition \ref{TFAE}, we have the following theorem :
\begin{thm}\label{isotropysubgroup}
Let $\Sigma$ be a graph of a function $F : \R^n \rightarrow \R$ with
$|\det DdF| = 1$. 
If the $\Aut(\Sigma)_0$ contains a transitive abelian
subgroup, then the followings hold.
\begin{enumerate}[\hspace{3mm}\rm (a)]
\item $\Aut(\Sigma)_0 = \A \rtimes \J$, where $\A$ is the
  simply transitive abelian subgroup and $\J$ is the isotropy subgroup
  of $\Aut(Sigma)_0$. 
\item $\J$ acts on the associated complete abelian LSA $\EuA$
  as an automorphism and similarity with respect to the Hessian type
  metric. 
\end{enumerate}
\end{thm}

A derivation $B$ on an LSA $\EuA$ is the linear map satisfying 
\[ B(ab) = (Ba)b + a(Bb), \]
for all $a, b \in \EuA$, that is, $B\lambda_a = \lambda_{Ba} +
\lambda_aB$ for all $a \in \EuA$.
Let $\Der(\EuA)$ be the set of all the derivation on $\EuA$, which is
the Lie algebra of the automorphism group of $\EuA$, that is,
$\Der(\EuA) = \Lie \Aut(\EuA)$.
 
\begin{definition}\label{def2}
For an LSA $\EuA$ with an inner product $H$, a map $B : \EuA
\rightarrow \EuA$ is called an \emph{infinitesimal similarity}  if it
satisfies the following :
\begin{eqnarray}\label{similarity}
& (Ba)'H b +  a'H B b  = \displaystyle s \: a'H b , \\
& (\mbox{ or } \langle Ba, b \rangle + \langle a, Bb \rangle =
s\langle a, b \rangle ) \nonumber
\end{eqnarray} 
where $s$ must be equal to $\displaystyle \frac{2}{n} \tr B$ for $n =
\dim \EuA$. 
\end{definition}

Let's denote the subalgebra of $\Der(\EuA)$ whose elements are
infinitesimal similarities  by $\sDer(\EuA)$. 
For an element $B \in \sDer(\EuA)$, if $\tr B = 0$, we call it
\emph{infinitesimal isometry derivation}.
We will denote by $\iDer(\EuA)$  the subalgebra of all
infinitesimal isometry derivations.
It readily follows that the codimension of $\iDer(\EuA)$ in $\sDer(\EuA)$ must be 0 or
\nolinebreak 1. 

\begin{proposition} \label{frakj}
Let $\Sigma$ be the graph of a function $F : \R^n \rightarrow \R$ with
$|\det DdF | = 1$, on which an abelian subgroup acts simply transitively.
Let $\frakj$ be the Lie algebra of the isotropy subgroup $\J$ of
$\Aut(\Sigma)_0$.
Then 
\[ \frakj = 
 \big\{  \begin{pmatrix} B   & 0 \\ 0 & s
\end{pmatrix} \; \big| \; B \in \sDer(\EuA) \mbox{ and } s =
\displaystyle \frac{2}{n} \tr B \big\}, \]
where $\EuA$ is the associated complete abelian LSA. 
Therefore the dimension of $\:\J$ is equal to the dimension of the Lie
algebra $\sDer(\EuA)$.
\end{proposition}

\begin{proof}
From the Lemma \ref{elem-isot}, the elements of $\frakj$ are given by
the matrix $\begin{pmatrix} B & 0 \\ 0 & s \end{pmatrix} \in \gl(n+1)$
where $s = \displaystyle \frac{2}{n} \tr B$. Then, from the
Proposition \ref{TFAE} (d), $B$ must be contained in $\sDer(\EuA)$. 
\end{proof}

If $\sDer(\EuA) \varsupsetneqq \iDer(\EuA)$, then the complete abelian LSA $\EuA$ has a
derivation $B$, which is an infinitesimal similarity as well with
$\displaystyle \tr B = \frac{n}{2}$, that is, $a'Hb = a'B'Hb + a'HBb$.
Hence the associated polynomial $F$ given in Theorem \ref{polynom} can be
represented as the following form, % :
\begin{eqnarray}
F(x) &=& x'H(\frac{1}{2}I - \frac{1}{3}\lambda_x + \frac{1}{4}(\lambda_x)^2
- \cdots)x \nonumber \\
&=& x'H(\frac{1}{2}x - \frac{1}{3}x^2 + \frac{1}{4}x^3 - \cdots) \nonumber \\
&=& x'B'H(\frac{1}{2}x - \frac{1}{3}x^2 + \frac{1}{4}x^3 - \cdots) +
x'H(\frac{1}{2}Bx - \frac{2}{3}xBx + \frac{3}{4}x^2Bx - \cdots) \nonumber \\ 
&=& x'H(Bx - xBx + x^2Bx - \cdots) \nonumber  \\
&=& x'H(I - \lambda_x + (\lambda_x)^2 - \cdots)Bx \nonumber  \\
&=& x' H (I + \lambda_x)^{-1} B x, \label{HBform}
\end{eqnarray}
by using the symmetry of $H$, the associativity  $(xy)'Hz = x'H(yz)$
and the derivation property, $Bx^k = kx^{k-1}Bx$ for
all $x, y, z \in \R^n$.

%%%%%%%%%%%%%%%%%%%%%%%%%%%%%%%%%%%%%%%%%%%%%%%%%%%%%%%%%%%%%%%%%%%%%
\section{Complete abelian LSA and Filiform LSA}\label{CompleteLSAandFili}

Let $\EuA$ be a complete abelian LSA 
and let $\EuA^i$'s be ideals given by the following :
\[ \EuA^1 = \EuA, \quad\quad  \EuA^{i+1} = \EuA \cdot \EuA^i \quad
\mbox{ for } i = 1,2,3, \cdots. \] 
Then there exists $m$ such that $\EuA^{m+1} = \{ 0 \}$ since
$\lambda_x$'s  are nilpotent for all $x \in \EuA$. 
In fact they make a
descending central series, that is,
\[ \EuA = \EuA^1 \supset \EuA^2 \supset \cdots \supset \EuA^m \supset
\EuA^{m+1} = \{ 0 \}. \] 
If $a \in \EuA^m$, then $b\cdot a \in \EuA^{m+1} = \{ 0 \}$ for any $b
\in \EuA$. This says that $\EuA^m$ is contained in the annihilator of
$\EuA$, $\Ann(\EuA) = \{ x \in \EuA | xy = 0 = yx, \mbox{ for all }y
\in \EuA \}$. 
Notice that $\Ann(\EuA)$ is a nontrivial ideal of a complete
abelian $\EuA$.

\begin{definition}[\cite{DO99}] Let $\EuA$ be a complete abelian
  LSA. A (vector space) basis of $\EuA$
\[ \{ e_{11}, \cdots, e_{1k_1}, e_{21},
\cdots, e_{2k_2}, \cdots, e_{m1}, \cdots, e_{mk_m} \} \]
is called \emph{adequate} if for all $i$, $e_{il}$ belongs to
$\EuA^i$ and $k_i + k_{i+1} + \cdots +k_m = \dim \EuA^i$.
\end{definition}

Let $\EuC_i$ be ideals defined by
\[ \EuC_1 = \Ann(\EuA), \quad\quad \EuC_{i+1} = \{ x \in \EuA \: | \:
yx  \in \EuC_i \mbox{ for any } y \in \EuA \}. \]
Then they make an
ascending central series, i.e.,
\[ \EuC_0 = \{ 0 \} \subset \EuC_1 \subset \EuC_2 \subset \cdots
. \]

\begin{lemma} \label{prop of center}
\begin{enumerate}[\hspace{3mm}\rm (a)]
\item $\EuA / \EuC_i$ is a complete abelian LSA for all $i = 1, 2,
  \cdots$.
\item $\Ann(\EuA / \EuC_i) = \EuC_{i+1} / \EuC_i$ for all $i = 1, 2,
  \cdots$.
\item There exists an integer $m'$ such that $\EuC_{m'} = \EuA$, so
  the ascending central series is terminated finitely.
\item $\EuC_{i} = \{ x \in \EuA \: | \:
 x \cdot \EuA^i = \{0\} \; \}$ for $i \in \{ 1, 2, \cdots, m' \}$.  
\end{enumerate}
\end{lemma}

\begin{proof}
(a) $\EuA / \EuC_i$ is also an abelian LSA since $\EuC_i$ is ideal.
Since $\lambda_x = \rho_x$ are nilpotent for all $x \in \EuA$,
$\lambda_{\bar{x}} = \rho_{\bar{x}}$ are also nilpotent for $\bar{x}
= x + \EuC_i \in \EuA / \EuC_i$. 
So we have $\tr \rho_{\bar{x}} = 0$ for all $x \in \EuA$ and $\EuA /
\EuC_i$ is complete.

(b) is clear from the definition. 

(c) From (a) and (b), $\EuC_{i+1} / \EuC_i$ are all nontrivial if
$\EuC_i \neq \EuA$. Since $\EuA$ is of finite dimensional, 
there exist $m'$ such that $\EuC_m' = \EuA $.

(d) Since
$x \in \EuC_i$ is equivalent to that $xy \in \EuC_{i-1}$ for all $y
\in \EuA$, it follows immediately by induction.
\end{proof}

\begin{proposition}\label{propofEuA}
\begin{enumerate}
\item[\rm(a)] $\EuA^{m-i} \subset \EuC_{i+1}$ for $i \in \{ 0, \cdots, m-1 \}$.

\item[\rm(b)] $m' =  m$.
\end{enumerate}
\end{proposition}

\begin{proof}
(a) Since $\EuA^{m-i} \EuA^{i+1}  \subset \EuA^{m+1} = \{ 0 \}$,
$\EuA^{m-i}$ is a subspace of $\EuC_{i+1}$ for all $i \in \{ 0, 1,2, \cdots, m-1
\}$ from Lemma \ref{prop of center} (d). 

(b)  From (a), $\EuA = \EuA^1 \subset \EuC_{m}$, hence we have $m' \leq m$.
If $m' < m$, then $\EuA^m = \EuA^{m-1}\EuA = \EuA^{m-1}\EuC_{m'} = \{
0 \}$. This is a contradiction.
\end{proof}

\begin{definition}
Let $V$ be a vector space and let $\langle, \rangle$ be an inner
product on $V$. For a subset $W$, the perpendicular
subspace $W^{\perp}$ is defined as the following :
\[ W^{\perp} = \{ x \in V \: | \: \langle x, y \rangle = 0, \mbox{ for
  all } y \in W \}. \]
\end{definition}

\begin{proposition}\label{generater}
Let $\EuA$ be a complete abelian LSA with a nondegenerate Hessian type
inner product $\langle, \rangle$.
Then $\EuC_i = (\EuA^{i+1})^{\perp}$ and the dimension of $\EuC_i /
\EuC_{i-1}$ is equal to the dimension 
of $\EuA^i / \EuA^{i+1}$ for $i = 1, 2, \cdots, m$. 
\end{proposition}

\begin{proof}
By the associativity (\ref{HessiantypeofAbel}) of the Hessian type
inner product, $\EuC_i \subset (\EuA^{i+1})^{\perp}$. Let $x \in
(\EuA^{i+1})^{\perp}$, then for any $a \in \EuA$ and $b \in \EuA^i$, we have
$\langle a, bx \rangle = \langle ab, x \rangle = 0.$ 
Since $\langle, \rangle$ is nondegenerate, $bx = 0$, that is, $x \in
\EuC_i$ from Lemma \ref{prop of center} (d).
Therefore $\EuC_i = (\EuA^{i+1})^{\perp}$. This says that $\dim
\EuC_i = \dim \EuA / \EuA^{i+1}$ for each $i = 1,2, \cdots,m$.
So we have $\dim \EuC_i / \EuC_{i-1} = \dim (\EuA / \EuA^{i+1}) /
(\EuA / \EuA^i) = \dim \EuA^{i} / \EuA^{i+1}$.
\end{proof}

Let $\EuH$ be an $n$-dimensional LSA. $\EuH$ is called a \emph{filiform
  LSA} if it satisfies :
\[ \dim \EuH^i = \dim \EuH_i = n+1-i\quad
\mbox{ for all } i \in \{ 1, \cdots, n \}, \]
where $\EuH^i = \EuH \EuH^{i-1}$ and $\EuH_i = \EuH_{i-1}\EuH$ with
$\EuH^1 = \EuH_1 = \EuH$.({\it
  cf.} \cite{DO99})
Note that $\EuH$ is a $n$-dimensional filiform LSA if and only if the
nilpotency of $\EuH$ is $n+1$, that is, $\EuH^{n+1} = \EuH_{n+1} = \{
0 \}$ and $\EuH^n \neq \{0 \}$, $\EuH_{n} \neq \{ 0 \}$.
It has the following properties \nolinebreak : 

\begin{proposition}[\cite{DO99}] \label{PropOfFili}
Let $\EuH$ be an $n$-dimensional filiform LSA.
\begin{enumerate}
\item[\rm(a)] The associated Lie algebra of $\EuH$ is nilpotent.
\item[\rm(b)] $\EuH$ is a complete LSA.
\item[\rm(c)] $\EuH^i$ is equal to $\EuH_i$ for all $i \in \{ 1, \cdots, n \}$.
\item[\rm(d)] $\EuH$ has a one-dimensional annihilator $\Ann(\EuH) = \EuH^n = \EuH_n$
\item[\rm(e)] For all $i \in \{1,\cdots, n\}$, $\EuH / \EuH^i$ is a filiform LSA.
\item[\rm(f)] There is an adequate basis $\{ e_1, \cdots, e_n \}$ satisfying
\[ e_1e_j = e_{j+1} \quad \mbox{ for all } j \in \{ 1,2,\cdots, n-1 \}. \]
\end{enumerate}
\end{proposition}

The adequate basis in the Proposition \ref{PropOfFili} (f) is
called \emph{strongly adequate}. Note that if $\EuA$ is an abelian
filiform LSA and $\{ e_1, \cdots, e_n \}$ is a strongly adequate
basis of $\EuA$, then $e_i = e_1^i$ for all $i \in \{ 1, \cdots, n
\}$ and we have $e_ie_j = e_{i+j}$.({\it cf.} \cite{DO99})

\begin{lemma} \label{fili-special}
Let $\EuA$ be an abelian filiform LSA with a strongly adequate
basis $\{ e_1, \cdots, e_n \}$, then we have the following {\rm :}
\begin{enumerate}
\item[\rm(a)] There is a nondegenerate Hessian type inner product
  $\langle ,\rangle$ on $\EuA$ which is given by
\begin{equation}\label{fili-metric}
\langle e_i, e_j\rangle = \begin{cases}  1, & i+j = n+1, \\ 0, & i+j
\neq n+1,
\end{cases}
\end{equation}
which is unique up to LSA automorphism, i.e., if $\langle ,\rangle$ is
a nondegenerate Hessian type inner product on $\EuA$, then there exists a
strongly adequate basis such that {\rm (\ref{fili-metric})} holds up to
sign.
\item[\rm (b)] The dimension of $\sDer(\EuA)$ is $1$ and $\iDer(\EuA) =
  \{ 0 \}$.
\end{enumerate}
\end{lemma}

\begin{proof}
(a) Let $\langle, \rangle$ be the inner product given in
(\ref{fili-metric}). Then $\langle, \rangle$ is obviously
nondegenerate. For any $i,j$ and $k \in \{ 1, \cdots, n\}$,
\begin{eqnarray*}
\langle e_i\cdot e_j , e_k \rangle = \langle e_{i+j}, e_k\rangle
= \begin{cases}  1, & i+j+k = n+1, \\ 0, & i+j+k \neq n+1,
\end{cases} \\
\langle e_i,  e_j\cdot e_k \rangle = \langle e_{i},
e_{j+k}\rangle = \begin{cases}  1, & i+j+k = n+1, \\ 0, & i+j+k \neq
n+1,\end{cases}
\end{eqnarray*}
hence we conclude that $\langle ,\rangle$ is of Hessian type.
Now, let $\langle, \rangle$ be a nondegenerate Hessian type inner product on
$\EuA$.
Then $\langle e_i, e_j \rangle = 0$ for $i+j > n+1$ since $\EuA^{n+1}
= \{ 0 \}$. 
Let $\langle e_i, e_j \rangle = s_k$ for $k = i+j \leq n+1$. 
Since $\langle, \rangle$ is nondegenerate, $s_{n+1}$ must not be zero.
By changing a basis inductively, we can obtain the
basis $\{ \bar{e}_1, \bar{e}_2, \cdots, \bar{e}_n \}$ satisfying
(\ref{fili-metric}) as follows :
First, by scaling $e_1$, we may assume that $s_{n+1} = \pm 1$. 
If $s_n \neq 0$, by putting $\bar{e_1} = e_1 \mp \displaystyle \frac{s_n}{n}
e_2$ and $\bar{e}_{j+1} = \bar{e}_1\bar{e}_j$ for $j \in \{ 1,2,
\cdots, n-1 \}$, we have that $\langle \bar{e}_i, \bar{e}_j \rangle =
0$ for $i +j = n$. 
If $s_n = \cdots = s_{k+1} = 0$ and $s_k \neq 0$, put 
$\bar{e}_1 = e_1 \mp \displaystyle \frac{s_k}{k}
e_{m-k+2}$ and $\bar{e}_{j+1} = \bar{e}_1\bar{e}_j$ for $j \in \{ 1,2,
\cdots, n-1 \}$. Then $\langle \bar{e}_i, \bar{e}_j \rangle =
0$ for $i +j = k, k+1, \cdots, n$.   

(b) Define a derivation $B_1 :\EuA \rightarrow \EuA$ by $B_1(e_1) =
e_1$, so that $B_1(e_j) = je_j$ for all $j = 2, 3, \cdots, n$.
Then $\tr B_1 = \displaystyle \frac{n(n+1)}{2}$ and
\[ \langle B_1(e_i), e_j\rangle + \langle e_i,
B_1(e_j) \rangle  
= \begin{cases}  n+1 , & i+j = n+1 \\ 0, & i+j \neq n+1
 \end{cases} \Bigg\}
= \frac{2}{n} \tr B_1 \langle e_i, e_j \rangle. \]
So $B_1$ is contained in $\sDer(\EuA)$, and hence $\dim \sDer(\EuA)
\geq 1$. %\\

For any $B \in \sDer(\EuA)$ and the strong adequate basis $\{ e_1,
\cdots, e_n \}$, $B(e_k) = B(e_1^k) = ke_{k-1}B(e_1)$, $(k =2, 3,
\cdots, n)$ by induction.  
Therefore if $B(e_1) = u_1e_1 + u_2e_2 + \cdots + u_ne_n$, then
$B(e_k) = ke_{k-1}B(e_1) = ke_{k-1} (u_1e_1 + u_2e_2 + 
\cdots + u_ne_n)$ for all $k = 2, 3, \cdots, n$ and hence $\tr B =
\displaystyle \frac{n(n+1)}{2} u_1$. 
Now assume that $0 \neq B \in \iDer(\EuA)$, then the $e_1$-coefficient of
$B(e_1)$ must be zero so that $B(e_1)$ is contained in $\EuA^2$.
In this case, there is a number $k$ such that $\langle B(e_1), e_k
\rangle \neq 0$ and $\langle B(e_1), e_{k+1} \rangle = 0$  since
$\langle, \rangle$ is nondegenerate and the basis is adequate. 
Note that $k$ is less than $n$ because $B(e_1) \in \EuA^2$. 
Since $B(e_k) = ke_{k-1}B(e_1)$, we have 
\[ 0 = \langle B(e_1), e_k \rangle + \langle e_1, B(e_k) \rangle =
(k+1) \langle B(e_1), e_k \rangle \neq 0. \]
This is a contradiction. Hence, $\iDer(\EuA) = \{ 0 \}$ and the
dimension of $\sDer(\EuA)$ has to be 1.
\end{proof}

The following key observation of this paper is made by Y. Choi\cite{Cho02}.

\begin{thm} \label{fili-Cayley}
Let $\EuA$ be a complete abelian LSA. Then $\EuA$ is  filiform if and only
if the associated polynomial $F$ of $\EuA$ with a nondegenerate
Hessian metric is the Cayley polynomial. 
\end{thm}

\begin{proof}
Let $\EuA$ be an abelian filiform LSA and $\{ e_1,
\cdots, e_n \}$ be a strongly adequate basis of $\EuA$. Then
from Lemma \ref{fili-special}, we have a Hessian type inner product $H$
and an infinitesimal similarity derivation $B$, which  are given
by the following matrices$\;$:
\[ H = \begin{pmatrix}
 &  &  &  & \; 1 \\
 &\;0\; &  &\; 1\; & \\
 &  &\bddots & & \\
 &\; 1\; &  &\;0\; &   \\
 1\; & &  &  &  \\
\end{pmatrix}, \quad
B = \begin{pmatrix}
\frac{1}{n+1} &  &  &  &  \\
 &\frac{2}{n+1}  &  & 0 & \\
 & &\ddots & &  \\
 & 0 &  &\frac{n-1}{n+1} &  \\
 &  &  &  & \frac{n}{n+1} \\
\end{pmatrix}. \]
The left multiplication $\lambda_x$ of $x = x_1e_1 + x_2e_2 + \cdots
+  x_ne_n$ and $(I + \lambda_x)^{-1}$ can be written as :
\[ \lambda_x = \begin{pmatrix}
0 &  &  &  &  \\
x_1 &0 &  &0 & \\
x_2 & x_1 &0 & & \\
\vdots  &\!\ddots &\!\ddots  &\!\ddots &   \\
x_{n-1} & \cdots& x_2 & x_1 & 0  \\
\end{pmatrix}, \quad
(I +\lambda_x)^{-1} = \begin{pmatrix}
1 &  &  &  &  \\
p_1 &1 &  &0 & \\
p_2 &p_1 &1 & & \\
\vdots  &\!\ddots &\!\ddots  &\!\ddots &   \\
p_{n-1} & \cdots& p_2 & p_1 & 1  \\
\end{pmatrix} \]
where $p_r = \displaystyle\sum_{i+j+\cdots + m = r} (-1)^{d}
\overbrace{x_ix_j \cdots x_m}^d $ for $r = 1,2,\cdots, n-1$. 
Then from (\ref{HBform}), the associated polynomial is given
by
\[ F(x) = x' H (I + \lambda_x)^{-1}B x 
= \sum_{d=2}^{n+1} (-1)^d \frac1d \sum_{i+j+\cdots
+ m = n+1} \overbrace{x_ix_j \cdots x_m}^d ,
\]
that is, $F$ is a Cayley polynomial.
Conversely, suppose that $\EuA$ is an $n$-dimensional complete
abelian LSA which is not filiform, then the nilpotency class of
$\lambda_x$ is less than $n$, that is, $\lambda_x^{n-1} = 0$ for all
$x \in \EuA$. 
This says that the degree of the associated
polynomial $F(x)$ given Theorem \ref{polynom} 
is less than $n+1$. Therefore the polynomial $F$ can not be
the Cayley polynomial.  
\end{proof}

Now we can restate the conjecture of Eastwood and Ezhov in terms of LSA

\begin{conjecture} [LSA version] \label{conjecture}
Let $\EuA$ be a complete abelian LSA with a nondegenerate Hessian type
inner product. If the dimension of $\sDer(\EuA)$ is $1$, then $\EuA$ is the
abelian filiform LSA.
\end{conjecture}

%%%%%%%%%%%%%%%%%%%%%%%%%%%%%%%%%%%%%%%%%%%%%%%%%%%%%%%%%%%%%%%%%%%%%%%%%%%%%
\section{Proof of the Main Theorem}

The affine $(n+1)$-space is divided into two connected components
by the hypersurface $\Sigma$ given by a function $F : \R^n \rightarrow
\R$. Choose the one which is lying above
$\Sigma$ and denote it by $\Omega$, then $\Sigma$ 
becomes a boundary of $\Omega$. Let $\Aut(\Omega)$ be the  subgroup of 
$\Aff(n+1)$ whose elements leave the domain $\Omega$ invariant. 

\begin{proposition}
$\Aut(\Omega) \subset \Aut(\Sigma)$ and $\Aut(\Sigma)_0 = \Aut(\Omega)_0$.
\end{proposition}

\begin{proof}
For any $ g \in \Aut(\Sigma)_0 \subset \Aff(n+1)$
and $x \in \Omega$,
if $gx \in \Sigma$ then put $y = gx$ and 
choose $g^{-1}y = z \in \Sigma$.
Then we have $gx = gz$, i.e., $g$ is not one-to-one. This contradicts that 
$g \in \Aff(n+1)$.
If $gx \in \Omega'$, the other connected component,
choose a path $g_t$ from $g_0 =e$ to $g_1=g$ in $\Aut(\Sigma)_0$.
Then $g_t x$ is a path from $x$ to $gx$ and 
we can choose $t_0$ such that $g_{t_0} x \in \Sigma$. 
This contradicts as above.
Hence we have $gx \in \Omega$, i.e., $g \in \Aut(\Omega)_0$.

Conversely, for any $g \in \Aut(\Omega) \subset \Aff(n+1)$
and $x \in \Sigma$,
if $gx \in \Omega$ then put $y = gx$ and  $g^{-1}y = x \in \Sigma$. 
This contradicts that $g \in \Aut(\Omega)$.
If $gx \in \Omega'$, choose a sequence $\{ x_n \} \subset \Omega$ 
which converges to $x$.
Then the sequence $\{ g\cdot x_n \} \subset \Omega$ converges to
$g\cdot x$,  hence there exists $x_N$ such that $gx_N \in \Omega'$. 
This is a contradiction.
Hence we have $gx \in \Sigma$, i.e., $g \in \Aut(\Sigma)$.
Similarly if $g \in \Aut(\Omega)_0$, then $g \in \Aut(\Sigma)_0$ by
considering a path $g_t$ with the same argument.
\end{proof}

The domain above the Cayley hypersurface is homogeneous, in fact, the
identity component of its automorphism group acts on the domain simply
transitively({\it cf.} Lemma \ref{fili-special} (b)). 
In the following, we will assume that the domain above the
hypersurface is also homogeneous. We call these hypersurfaces
\emph{extensible homogeneous}. 

\begin{lemma} \label{isotropy}
Let $\Sigma$ be a graph of a function $F : \R^n \rightarrow \R$ with
$|\det DdF | =1$, on which the
abelian subgroup $\A$ acts simply transitively. Let $\J$ be the
isotropy subgroup of $\Aut(\Sigma)$ and let $\EuA$ be the associated
complete abelian LSA with Hessian type inner product $H$. 
Then $\Sigma$ is extensible homogeneous if and only if the action of
$\J$ on $(\EuA, H)$ contains a similarity which is not an isometry. 
\end{lemma}

\begin{proof}
If the action of $\J$ contains only isometries, then from
Proposition \ref{TFAE} any element $j \in \J$ is represented by $j =
\begin{pmatrix} A & 0 \\ 0 & 1 \end{pmatrix}$,
hence $j\cdot e_{n+1} = e_{n+1}$ where $e_{n+1} = (0, \cdots, 0, 1)$. 
This says that $\Aut(\Sigma)_0 \cdot e_{n+1} \subset \Sigma_1 = \{ (x,
F(x) + 1 ) \: | \: (x, F(x)) \in \Sigma \}$, 
which contradicts that $\Aut(\Sigma)_0$ acts transitively on 
$\Omega$.

Conversely, if the action of $\J$ on $(\EuA, H)$ contains a
similarity which is not an isometry, then $\J$ contains 1-parameter
subgroup such that $\J \cdot e_{n+1} = \{ (0,
\cdots, 0, t) \: | \: t \in \R_+ \}$ from Proposition \ref{frakj}. 
Moreover $\A$ acts on $\Sigma_t
= \{ (x, F(x) + t) \: | \: (x, F(x)) \in \Sigma, \; t \in \R_+ \}$
simply transitively by using the matrix form in Lemma \ref{gpelementofAbel}.
Therefore $\Aut(\Sigma)_0 = \A \rtimes \J$ acts on $\Omega$ transitively.
\end{proof}

From Lemma \ref{isotropy}, a complete abelian LSA $\EuA$, which is associated
to the extensible homogeneous hypersurface, has a derivation $B$ such
that $B \in \sDer(\EuA)$ and $B \not\in \iDer(\EuA)$, i.e.,  $\tr B \neq 0$.
In this case, for a proof of the Conjecture, it is enough to show
that if $\EuA$ is not filiform, $\iDer(\EuA)$ is not trivial. Note
that,  if $\EuA$ is an abelian filiform LSA, then the dimension of
$\sDer(\EuA)$ is 1 and $\iDer(\EuA) = \{ 0 \}$ from Lemma \ref{fili-special}.

\begin{proposition}\label{CenterOfFili}
Let $\EuA$ be a complete abelian LSA with a nondegenerate Hessian
type inner product $\langle ,\rangle$, then $\EuA$ is filiform  if and
only if $\dim \Ann(\EuA)= 1$.
\end{proposition}

\begin{proof}
If $\EuA$ is filiform, then $\dim \Ann(\EuA)=1$ from the Proposition
\ref{PropOfFili}. Now suppose that $\Ann(\EuA)= \spann\{\, x \} = \EuA^m$.
Choose $e_1 \in \EuA$ such that $\langle e_1, x \rangle = 1$.
Because $\EuA^2$ is perpendicular to $\EuC_1 = \Ann(\EuA)$, we have $e_1
\not\in \EuA^2$.
Since $\dim \EuA / \EuA^2 = \dim \EuC_1 = 1$ from Lemma
\ref{generater}, $\EuA$ can be represented as the following
\nolinebreak :
\[ \EuA = V + \EuA^2, \]
where $V$ is the 1-dimensional vector space spanned by $\{ e_1 \}$.
Then, inductively, we have $\EuA^i = V^i + \EuA^{i+1}$ for all $i$, and
hence
\[ \EuA = V + V^2 + \cdots + V^m, \]
where $V^i = \spann \{ e_1^i \}$ for $1 \leq i \leq m$. 
Observe that $\{ e_1, e_1^2, \cdots, e_1^m \}$ is linearly independent 
and hence becomes a basis of $\EuA$, that is, $\EuA$ is a filiform
LSA. 
\end{proof}

\begin{lemma} \label{DimOfCenter}
Let $\EuA$ be a complete abelian LSA with a nondegenerate Hessian
type inner product $\langle ,\rangle$,  and let $d = \dim \Ann(\EuA)$
be greater than or equal to $2$. 
Then the Lie algebra $\iDer(\EuA)$ is nontrivial.
\end{lemma}

\begin{proof}
Let $\EuZ = \spann \{ e_{n-1}, e_n \}$ be a
subspace of $\Ann(\EuA)$.

{\bf (Case 1)} $\EuZ \subset \EuZ^{\perp}$. Consider a linear map
$\psi : \EuA \rightarrow \R^2, \; x \mapsto (\langle x,
e_{n-1}\rangle, \langle x, e_n \rangle)$. 
Then $\ker \psi = \EuZ^{\perp}$. 
Since $\EuZ^{\perp}$ has codimension 2, $\phi$ is onto and hence there
exist elements $e_1, e_2 \in \EuA$ 
such that 
\begin{eqnarray*}
& \langle e_1, e_{n-1} \rangle = 0, \quad \langle e_1, e_n \rangle = 1
\\
& \langle e_2, e_{n-1} \rangle = 1, \quad \langle e_2, e_n \rangle = 0 
\end{eqnarray*} 
Denote $V = \spann \{ e_1, e_2 \}$, then $\EuA = V + \EuZ^{\perp}$.
Note that $\EuZ \subset \Ann(\EuA)$ implies $\EuZ^{\perp} \supset
\Ann(\EuA)^{\perp} = \EuA^2$.
Now define a map $B : \EuA \rightarrow \EuA$ by 
\[ B(e_1) = - e_{n-1}, \; B(e_2) = e_n, \; B(x) = 0 \mbox{ for } x \in
\EuZ^{\perp}. \]
Then $B$ is a derivation since $B(\EuA) = \EuZ \subset \Ann(\EuA)$ and
$B(\EuA^2) = \{ 0 \}$.
Let $\{ e_1, e_2, e_3, \cdots, e_{n-2}, e_{n-1}, e_n \}$ be a basis of
$\EuA$ where $ e_3, \cdots, e_{n-2} \in \EuZ^{\perp}$. Then 
the matrix form $H$ of the Hessian type inner product $\langle,
\rangle$ and the derivation $B$ are given as follows,
\[ H = \begin{pmatrix}
 &  & & & 0 & 1 \\
 &  *&  &* & 1 & 0\\
 &  &  & & 0 & 0\\
 &  *&  & * & \vdots & \vdots \\
0 & 1 & 0 &\cdots& 0 & 0 \\
1 & 0 & 0 &\cdots& 0 & 0 \\
\end{pmatrix}, \quad  B = \begin{pmatrix}
0 & 0 & 0 & \cdots &0 & 0 \\
0 & 0 & 0 & \cdots &0 & 0\\
\vdots &\vdots &\vdots &\vdots &\vdots &\vdots \\
0 & 0 & 0 & \cdots &0 & 0\\
-1 & 0 & 0 & \cdots &0 & 0 \\
0 & 1 & 0 & \cdots &0 & 0 \\
\end{pmatrix}. \] 
By matrix calculation, we have $B^t H + H B = 0$.
So $B$ is contained in $\iDer(\EuA)$, that is, $\iDer(\EuA) \neq \{0 \}$.

{\bf (Case 2)} $\EuZ \not\subset \EuZ^{\perp}$ and $\EuZ \cap
\EuZ^{\perp} \neq \{ 0 \}$. 
By coordinate change, we may assume that
$e_{n-1} \in \EuZ \cap \EuZ^{\perp}$ and $\langle e_n , e_n \rangle =
1$ using  $\EuZ \not\subset \EuZ^{\perp}$.
Put $\EuZ_1 = \spann \{ e_n \}$. 
Then since $\EuZ_1$ is an ideal, by the associativity of Hessian type
inner product, $\bar{\EuA} = \EuZ_1^{\perp}$ becomes
an ideal of $\EuA$. 
Since $\EuZ_1$ is a nondegenerate subspace,  $\EuA = \bar{\EuA} \oplus
\EuZ_1$. 
Note
that $\EuZ \cap \EuZ^{\perp} \subset \bar{\EuA}$. Choose $e_1 \in
\bar{\EuA}$ such that $\langle e_1, e_{n-1} \rangle = 1$. 
In this case define a derivation $B$ by 
\[ B(e_1) = -e_n, \; B(e_n) = e_{n-1}, B(x) = 0 \mbox{ for } x \in
\EuZ^{\perp}. \]
Let $\{ e_1, e_2, e_3, \cdots, e_{n-2}, e_{n-1}, e_n \}$ be a basis of
$\EuA$ where $ e_2, e_3 \cdots, e_{n-2} \in \EuZ^{\perp}$. Then 
the matrix forms $H$ of the Hessian type inner product $\langle,
\rangle$ and the derivation $B$ are given as follows
\[ H = \begin{pmatrix}
 &  & & & 1 & 0 \\
 &  *&  &* & 0 & 0\\
 &  &  & & 0 & 0\\
 &  *&  & * & \vdots & \vdots \\
1 & 0 & 0 &\cdots& 0 & 0 \\
0 & 0 & 0 &\cdots& 0 & 1 \\
\end{pmatrix}, \quad  B = \begin{pmatrix}
0 & 0 & 0 & \cdots &0 & 0 \\
0 & 0 & 0 & \cdots &0 & 0\\
\vdots &\vdots &\vdots &\vdots &\vdots &\vdots \\
0 & 0 & 0 & \cdots &0 & 0\\
0 & 0 & 0 & \cdots &0 & 1 \\
-1 & 0 & 0 & \cdots &0 & 0 \\
\end{pmatrix}. \] 
By matrix calculation, we have $B^t H + H B = 0$.
So $B$ is contained in $\iDer(\EuA)$, that is, $\iDer(\EuA) \neq \{0 \}$.

{\bf (Case 3)} $\EuZ \cap \EuZ^{\perp} = \{ 0 \}$. 
In this case, the restriction of the inner
product on $\EuZ$ is nondegenerate. 
So $\EuZ^{\perp}$ becomes a nondegenerate ideal
again by the associativity of $\langle, \rangle$.
Therefore we have an algebra direct sum $\EuA = \EuZ^{\perp} \oplus \EuZ$. 
Since $\EuZ$ is a trivial LSA, $\iDer(\EuZ) = \frako(2)$ or
$\frako(1,1)$, which are clearly nontrivial. 
For any $\bar{B} \in \iDer{\EuZ}$ the map $B : \EuA \rightarrow \EuA$
given by 
\[ B(x) = \begin{cases} \bar{B}(x) & x \in \EuZ, \\ 0 & x \in
  \EuZ^{\perp} \end{cases} \]
is contained in $\iDer(\EuA)$. This says that $\iDer(\EuA) \neq \{0 \}$.
\end{proof}

\begin{corollary} \label{equivdim}
Let $\EuA$ be a complete abelian LSA which is associated to the
nondegenerate extensible homogeneous hypersurface.
Then the followings are equivalent $:$
\begin{enumerate}[\hspace{4mm} \rm (a)]
 \item $\EuA$ is filiform.
 \item $\dim \Ann(\EuA)= 1$.
 \item $\dim \sDer(\EuA) = 1$.
\end{enumerate}
\end{corollary}

\begin{proof}
Since this complete abelian LSA $\EuA$ has a nondegenerate Hessian
type inner product, Proposition \ref{CenterOfFili} shows that $ \EuA$
is filiform if and only if $\dim \Ann(\EuA) = 1$, and the fact that
$\EuA$ is filiform implies that $\dim \sDer(\EuA) = 1 $ from Lemma
\ref{fili-special}. 
Because the LSA $\EuA$ is associated to the extensible hypersurface,
Lemma \ref{isotropy} says that $\dim \sDer(\EuA) / \iDer(\EuA) = 1$. 
Finally, Lemma \ref{DimOfCenter} shows that if $\dim \sDer(\EuA) =1$,
then $\dim \Ann(\EuA) = 1$.
\end{proof}

\begin{mainthm} %\label{main}
Let $\Sigma$ be a nondegenerate hypersurface given by a function on
$\R^n$. 
Then $\Sigma$ is a Cayley hypersurface if and only if the followings
are satisfied \nolinebreak$:$
\begin{enumerate}[\hspace{4mm} \rm (E1)]
\item $\Sigma$ admits a transitive abelian group $\A$ of affine
motions.
\item $\Aut(\Sigma)$  has a 1-dimensional isotropy group.
\item Affine normals to $\Sigma$ are everywhere parallel.
\item The domain above the hypersurface is also homogeneous.
\end{enumerate}
\end{mainthm}

\begin{proof}
For an $n$-dimensional abelian filiform LSA $\EuA$ with Hessian type
inner product $H=\langle, \rangle$ given by (\ref{fili-metric}),  let
$\A$ be the associated abelian Lie subgroup of $\Aff(n+1)$ obtained as
in Theorem \ref{one-to-one}({\it see} the discussion following Remark
\ref{remark}), that is, $\Lie \A$ is the associated Lie algebra of
$\EuA$.  
Then Theorem \ref{fili-Cayley} shows that $\A$ acts on the Cayley
hypersurface $\Sigma$ simply transitively. 
The conditions, (E2) and (E4) follow from Lemma \ref{fili-special}(b) and
Lemma \ref{isotropy}. 
From Theorem \ref{one-to-one}, the affine normals of $\Sigma$ are
parallel. 

Conversely, the transitive abelian group $\A$ acts in fact simply
transitively on $\Sigma$ from Proposition \ref{simply}. 
Since $\Sigma$ is a global graph whose affine normals are parallel,
$\A$ gives us a complete abelian LSA $\EuA$ with nondegenerate Hessian
type inner product by Theorem \ref{one-to-one}. 
From {\rm (E2)} with Proposition \ref{frakj}, $\dim \sDer(\EuA) = 1$ and
hence by Corollary \ref{equivdim}, $\EuA$ is filiform.  
Then by Theorem \ref{fili-Cayley}, the polynomial $F$ obtained from 
$\EuA$ is the Cayley polynomial, that is, $\Sigma$ is a Cayley
hypersurface.
\end{proof}

It is unclear for us at this moment that (E1) and (E2) implies (E4) so that
the condition (E4) can be  eliminated in the Main Theorem.  
In this paper we consider a hypersurface given by a graph of a global
function defined on {\it all} of $\R^n$. 
It would be very interesting to know whether only the conditions
(E1), (E2) and (E3) would imply the globalness of the function.
The study of these two problems will give a complete answer to the
conjecture of Eastwood and Ezhov.
Furthermore the classification problem of affine homogeneous
hypersurfaces could be further studied through the LSA
with Hessian type inner product as suggested in this paper.

\end{document}